\documentclass[11pt]{article}
\usepackage{graphicx,array,amssymb,psfrag,fullpage}

\usepackage{algorithmic,algorithm}

\newcommand{\BEAS}{\begin{eqnarray*}}
\newcommand{\EEAS}{\end{eqnarray*}}
\newcommand{\BEA}{\begin{eqnarray}}
\newcommand{\EEA}{\end{eqnarray}}
\newcommand{\BEQ}{\begin{equation}}
\newcommand{\EEQ}{\end{equation}}
\newcommand{\BIT}{\begin{itemize}}
\newcommand{\EIT}{\end{itemize}}
\newcommand{\BNUM}{\begin{enumerate}}
\newcommand{\ENUM}{\end{enumerate}}

\newcommand{\BA}{\begin{array}}
\newcommand{\EA}{\end{array}}

\newcommand{\ones}{\mathbf 1}

\newcommand{\reals}{{\mbox{\bf R}}}

\newcommand{\symm}{{\mbox{\bf S}}}  

\newcommand{\Rank}{\mathop{\bf Rank}}
\newcommand{\NumRank}{\mathop{\bf NumRank}}
\newcommand{\Card}{\mathop{\bf Card}}
\newcommand{\Tr}{\mathop{\bf Tr}}
\newcommand{\diag}{\mathop{\bf diag}}
\newcommand{\lambdamax}{{\lambda_{\rm max}}}

\newcommand{\Expect}{\textstyle\mathop{\bf E{}}}
\newcommand{\Prob}{\mathop{\bf Prob}}

\newcommand{\QED}{~~\rule[-1pt]{6pt}{6pt}}

\newcommand{\argmin}{\mathop{\rm argmin}}

\newcommand{\vect}{\mathop{\bf vec}}

\newtheorem{theorem}{Theorem}

\newtheorem{remark}[theorem]{Remark}

\newcounter{exno}

\newenvironment{proof}{\textbf{Proof.}}{\QED\bigskip}

{\begin{quote}}{\end{quote}}

\makeatletter
\long\def\@makecaption#1#2{
   \vskip 9pt 
   \begin{small}
   \setbox\@tempboxa\hbox{{\bf #1:} #2}
   \ifdim \wd\@tempboxa > 5.5in
        \begin{center}
        \begin{minipage}[t]{5.5in}
        \addtolength{\baselineskip}{-0.95pt}
        {\bf #1:} #2 \par
        \addtolength{\baselineskip}{0.95pt}
        \end{minipage}
        \end{center}
   \else 
    \hbox to\hsize{\hfil\box\@tempboxa\hfil}  
   \fi
   \end{small}\par
}
\makeatother

\newcounter{oursection}

\newcounter{lecture}

\newtheorem{lemma}{Lemma}
\begin{document}

\title{Subsampling Algorithms for Semidefinite Programming.}

\author{Alexandre d'Aspremont\thanks{CMAP, \'Ecole Polytechnique, UMR CNRS 7641, Palaiseau, France. \texttt{alexandre.daspremont@m4x.org}}}

\maketitle

\begin{abstract}
We derive a stochastic gradient algorithm for semidefinite optimization using randomization techniques. The algorithm uses subsampling to reduce the computational cost of each iteration and the subsampling ratio explicitly controls  granularity, i.e. the tradeoff between cost per iteration and total number of iterations. Furthermore, the total computational cost is directly proportional to the complexity (i.e. rank) of the solution. We study numerical performance on some large-scale problems arising in statistical learning.
\end{abstract}


\section{Introduction}
Beyond classic combinatorial relaxations \cite{Goem95}, semidefinite programming has recently found a new stream of applications in machine learning \cite{Lanc02}, geometry \cite{Wein06}, statistics \cite{dAsp06b} or graph theory \cite{Sun05}. All these problems have a common characteristic: they have relatively low precision targets but form very large semidefinite programs for which obtaining second order models is numerically hopeless which means that Newton based interior point solvers typically fail before completing even a single iteration. Early efforts focused on exploiting structural properties of the problem (sparsity, block patterns, etc), but this has proven particularly hard for semidefinite programs. For very large problem instances, first-order methods remain at this point the only credible alternative. This follows a more general trend in optimization which seeks to significantly reduce the \emph{granularity} of solvers, i.e. reduce the per iteration complexity of optimization algorithms rather than their total computational cost, thus allowing at least some progress to be made on problems that are beyond the reach of current algorithms. 

In this work, we focus on the following spectral norm minimization problem
\BEQ\label{eq:min-maxeig}
\BA{ll}
\mbox{minimize} & \left\|\sum_{j=1}^p y_j A_j +C\right\|_2-b^Ty\\
\mbox{subject to} & y \in Q,
\EA\EEQ
in the variable $y\in \reals^p$, with parameters $A_j\in\symm_n$, for $j=1,\ldots,p$, $b\in\reals^p$ and $C\in \symm_n$, where $Q$ is a compact convex set. Throughout the paper, we also implicitly assume that the set $Q\subset\reals^p$ is simple enough so that the complexity of projecting $y$ on $Q$ is relatively low compared to the other steps in the algorithm.

The idea behind this paper stems from a recent result by \cite{Judi07}, who used a mirror descent stochastic approximation algorithm for solving bilinear matrix games (see \cite{Nest09}, \cite{Poly92} or \cite{Nemi83} for more background), where subsampling is used to perform matrix vector products and produce an approximate gradient. Strikingly, the algorithm has a total complexity of $O(n\log n/\epsilon^2)$, when the problem matrix is $n \times n$, hence only requires access to a negligible proportion of the matrix coefficients as the dimension $n$ tends to infinity. A similar subsampling argument was used in \cite{Judi08a} to solve a variational inequality representation of maximum eigenvalue minimization problems. 

In parallel, recent advances in large deviations and random matrix theory have produced a stream of new randomization results for high dimensional linear algebra (see \cite{Frie04,Drin06,Achl07, Vemp09} among many others), motivated by the need to perform these operations on very large scale, sometimes streaming, data sets in applications such as machine learning, signal processing, etc. Similar subsampling techniques have been successfully applied to support vector machine classification \cite{Kuma08} or Fourier decomposition. Randomization results were used in \cite{Aror07} to produce complexity bounds for certain semidefinite programs arising in combinatorial relaxations of graph problems. Randomization was also used in \cite{Burk02} and \cite{Burk05} to approximate subdifferentials of functions that are only differentiable almost everywhere. A randomized algorithm for semidefinite programming based on random walk techniques was also developed in \cite{Shch07}. Finally, a recent stochastic version by \cite{Lan09} of the algorithm in \cite{Nest07} has the potential to improve the complexity bounds provided by the method in \cite{Judi07}.

Our contribution here is to further reduce the granularity of first-order semidefinite programming solvers by combining subsampling procedures with stochastic approximation algorithms to derive stochastic gradient methods for spectral norm minimization with very low complexity per iteration. In practice, significantly larger per iteration complexity and memory requirements mean that interior point techniques often fail to complete a single iteration on very large problem instances. CPU clock also runs much faster than RAM, so operations small enough to be performed entirely in cache (which runs at full speed) are much faster than those requiring larger data sets. Solver performance on very large problem instances is then often more constrained by memory bandwidth than clock speed, hence everything else being equal, algorithms running many cheap iterations will be much faster than those requiring fewer, more complex ones. Here, subsampling techniques allow us to produce semidefinite optimization algorithms with very low cost per iteration, where all remaining $O(n^2)$ operations have a small constant and can be performed in a single pass over the data. 

We also observe that the relative approximation error in computing the spectral norm (or trace norm) of a matrix using subsampling is directly proportional to the numerical rank of that matrix, hence another important consequence of using subsampling techniques to solve large-scale semidefinite programs is that the total complexity of running the algorithm becomes explicitly dependent on the complexity (i.e. rank) of its solution. Most algorithms exploiting the fact that the solution has low rank are not convex (e.g. alternating minimization).

The paper is organized as follows. Section~\ref{s:random} surveys some key results on randomized linear algebra and spectral norm approximations. In Section~\ref{s:stoch-opt} we then derive a stochastic approximation algorithm for spectral norm minimization with very low cost per iteration and discuss some extensions to statistical learning problems. Finally, we present some numerical experiments in Section~\ref{s:numexp}.

\subsection*{Notation}
We write $\symm_n$ the set of symmetric matrices of dimension $n$. For a matrix $X\in\reals^{m\times n}$, we write $\|X\|_F$ its Frobenius norm, $\|X\|_{\mathrm{tr}}$ its trace norm, $\|X\|_2$ its spectral norm, $\sigma_i(X)$ its $i$-th largest singular value and let $\|X\|_\infty=\max_{ij}|X_{ij}|$, while $X^{(i)}$ is the $i$-th column of the matrix $X$ and $X_{(i)}$ its $i$-th row. We write $\vect(X)$ the vector of $\reals^{mn}$ obtained by stacking up the columns of the matrix $X$ and $\NumRank(X)$ the numerical rank of the matrix $X$, where $\NumRank(X)=\|X\|_F^2/\|X\|_2^2$. Finally, when $x\in\reals^n$ is a vector, we write~$\|x\|_2$ its Euclidean norm, while $\|\cdot\|$ is a general norm on $\reals^m$ and $\|\cdot\|_*$ its dual norm.

\section{Randomized linear algebra}
\label{s:random}
In this section, we survey several results by \cite{Drin06} which, after a single pass on the data, sample columns to approximate matrix products and produce low rank matrix approximations with a complexity of $O(sn)$ where $s$ is the sampling rate.

\subsection{Randomized matrix multiplication}

\begin{algorithm} 
\caption{Matrix multiplication} 
\label{alg:matrix-mult} 
\begin{algorithmic} [1]
\REQUIRE $A\in\reals^{m \times n}$, $B\in\reals^{n \times p}$ and $s$ such that $1 \leq s \leq n$.
\STATE Define a probability vector $q\in\reals^n$ such that 
\[
q_i=\frac{\|A^{(i)}\|_2\|B_{(i)}\|_2}{\sum_{j=1}^n \|A^{(j)}\|_2\|B_{(j)}\|_2}, \quad i=1,\ldots,n.
\]
\STATE Define subsampled matrices $C\in\reals^{m \times s}$ and $R\in\reals^{s \times p}$ as follows.
\FOR{$i=1$ to $s$} 
\STATE Pick $j\in[1,n]$ with $\mathbf{P}(j=l)=q_l$.
\STATE Set $C^{(i)}=A^{(j)}/\sqrt{sq_{j}}$ and $R_{(i)}=B_{(j)}/\sqrt{sq_{j}}$.
\ENDFOR 
\ENSURE Matrix product $CR$ approximating $AB$.
\end{algorithmic} 
\end{algorithm} 

By construction, $\Expect[CR]=AB$, and the following randomization result from \cite{Drin07} controls the precision of the approximations in algorithm \ref{alg:matrix-mult}. 

\begin{lemma}\label{ref:lem-col-sample}
Let $A\in\reals^{m \times n}$, $B\in\reals^{n \times p}$, given a subsampling rate $s$ such that $1 \leq s \leq n$, suppose that $C\in\reals^{m \times s}$ and $R\in\reals^{s \times p}$ are computed according to algorithm \ref{alg:matrix-mult} above, then
\[
\Expect[\|AB-CR\|_F^2]\leq \frac{1}{s}{\|A\|_F^2\|B\|_F^2}
\]
and if $\beta\in[0,1]$ with $\eta=1+\sqrt{8\log(1/\beta)}$ then
\[
\|AB-CR\|_F^2\leq \frac{\eta^2}{s}{\|A\|_F^2\|B\|_F^2}
\]
with probability at least $1-\beta$.
\end{lemma}
\begin{proof}
See Theorem 1 in \cite{Drin07}.
\end{proof}

Note that using the adaptive probabilities $q_i$ is crucial here. The error bounds increase by a factor $n$ when $q_i=1/n$ for example.

\subsection{Randomized low-rank approximation}

\begin{algorithm} 
\caption{Low-rank approximation} 
\label{alg:low-rank} 
\begin{algorithmic} [1]
\REQUIRE $X\in\reals^{m \times n}$ and $k,s$ such that $1\leq k \leq s < n$.
\STATE Define a probability vector $q\in\reals^n$ such that $q_i={\|X^{(i)}\|^2_2}/{\|X\|_F^2}$, for $i=1,\ldots,n$.
\STATE Define a subsampled matrix $S\in\reals^{m \times s}$ as follows.
\FOR{$i=1$ to $s$} 
\STATE Pick an index $j\in[1,n]$ with $\mathbf{P}(j=l)=q_l$.
\STATE Set $S^{(i)}=X^{(j)}/\sqrt{sq_{j}}$.
\ENDFOR 
\STATE Form the eigenvalue decomposition $S^TS=Y\diag(\sigma)Y^T$ where $Y\in\reals^{s \times s}$ and $\sigma \in \reals^s$.
\STATE Form a matrix $H\in\reals^{m\times k}$ with $H^{(i)}=SY^{(i)}/\sigma_i^{1/2}$.
\ENSURE Approximate singular vectors $H^{(i)}$, $i=1,\ldots,k$.
\end{algorithmic} 
\end{algorithm} 

Algorithm \ref{alg:low-rank} below computes the leading singular vectors of a smaller matrix $S$, which is a subsampled and rescaled version of $X$. Here, the computational savings come from the fact that we only need to compute singular values of a matrix of dimension~$m\times s$ with $s\leq n$. Recall that computing $k$ leading eigenvectors of a symmetric matrix of dimension $s$ only requires matrix vector products, hence can be performed in $O(ks^2\log s)$ operations using iterative algorithms such as the power method or Lanczos method. The complexity is in fact $O(kms\log m)$ in our case because an explicit factorization of the matrix is known, (see the appendix for details, as usual we omit the precision target in linear algebra operations, implicitly assuming that it is much finer than $\epsilon$), so that the cost of computing $k$ leading singular vectors of a matrix of size~$m\times s$ is $O(ksm\log m)$.

This means that, given the probabilities $q_i$, the total cost of obtaining $k$ approximate singular vectors using algorithm~\ref{alg:low-rank} is $O(ksm\log m)$ instead of $O(knm \log m)$ for exact singular vectors. Of course, computing $q_i$ requires $mn$ operations, but can be done very efficiently in a single pass over the data. We now recall the following result from \cite{Drin06} which controls the precision of the approximations in algorithm \ref{alg:low-rank}.

\begin{lemma}\label{ref:lem-col-vec}
Let $X\in\reals^{m \times n}$ and $1\leq k \leq s < n$. Given a precision target $\epsilon>0$, if $s\geq 4/\epsilon^2$ and $H\in\reals^{m \times k}$ is computed as in algorithm \ref{alg:low-rank}, we have
\[
\Expect[\|X-H_kH_k^TX\|_2^2]\leq \|X-X_k\|_2^2 + \epsilon \|X\|_F^2
\]
and if in addition $s>4\eta^2/\epsilon^2$ where $\eta=1+\sqrt{8\log(1/\beta)}$ for $\beta\in[0,1]$, then
\[
\|X-H_kH_k^TX\|_2^2 \leq \|X-X_k\|_2^2 + \epsilon \|X\|_F^2
\]
with probability at least $1-\beta$, where $X_k$ is the best rank $k$ approximation of $X$. 
\end{lemma}
\begin{proof}
See Theorem 4 in \cite{Drin06}.
\end{proof}

An identical precision bound holds in the Frobenius norm when $s\geq 4k/\epsilon^2$. We now adapt these results to our setting in the following lemma, which shows how to approximate the spectral radius of a symmetric matrix $X$ using algorithm \ref{alg:low-rank}.

\begin{lemma}\label{lem:col-approx}
Let $X\in\reals^{m \times n}$ and $\beta\in[0,1]$. Given a precision target $\epsilon>0$, construct a matrix $S\in\reals^{m \times s}$ by subsampling the columns of $X$ as in algorithm \ref{alg:low-rank}. Let $\eta=1+\sqrt{8\log(1/\beta)}$ and
\BEQ\label{eq:col-samp-rate}
s=\eta^2 \frac{\|X\|_2^2}{\epsilon^2 }\NumRank(X)^2
\EEQ
we have
\[
\Expect[|\|S\|_2-\|X\|_2|]\leq \epsilon
\]
and
\[
|\|S\|_2-\|X\|_2|\leq \epsilon
\]
with probability at least $1-\beta$.
\end{lemma}
\begin{proof}
Using the Hoffman-Wielandt inequality (see \cite[Th. 3.1]{Stew90} or the proof of \cite[Th.2]{Drin06} for example) we get
\[
|\|S\|_2^2-\|X\|_2^2|\leq \|SS^T-XX^T\|_F
\]
hence 
\[
|\|S\|_2-\|X\|_2|\leq \|SS^T-XX^T\|_F/\|X\|_2
\]
and Jensen's inequality together with the matrix multiplication result in Lemma \ref{ref:lem-col-sample} yields
\[
\Expect[\|SS^T-XX^T\|_F] \leq \displaystyle \frac{\|X\|_F^2}{\sqrt{s}}
\]
and
\[
\|SS^T-XX^T\|_F \leq \frac{\eta \|X\|_F^2}{\sqrt{s}}
\]
with probability at least $1-\beta$. Combining these two inequalities with the sampling rate in~(\ref{eq:col-samp-rate}) 
\[
s=\eta^2 \frac{\|X\|_F^4}{\epsilon^2 \|X\|_2^2}
\]
yields the desired result.
\end{proof}

The subsampling rate required to achieve a precision target $\epsilon$ has a natural interpretation. Indeed
\[
s=\eta^2 \frac{\|X\|_2^2}{\epsilon^2 }\NumRank(X)^2
\]
is simply the squared ratio of the numerical rank of the matrix $X$ over the relative precision target~${\epsilon}/{\|X\|_2}$, times a factor $\eta^2$ controlling the confidence level. The numerical rank $\NumRank(X)$ always satisfies $1 \leq\NumRank(X)={\|X\|_F^2}/{\|X\|_2^2}  \leq \Rank(X)$ and can be seen as a stable relaxation of the rank of the matrix $X$ (see \cite{Rude07} for a discussion). Note also that, by construction, the subsampled matrix always has lower rank than the matrix $X$. The expectation bound is still valid if we drop the factor $\eta$ in~(\ref{eq:col-samp-rate}).

\section{Stochastic approximation algorithm} 
\label{s:stoch-opt} Below, we will use a stochastic approximation algorithm to solve problem (\ref{eq:min-maxeig}) when the gradient is approximated using the subsampling algorithms detailed above. We focus on a stochastic approximation of problem (\ref{eq:min-maxeig}) written
\BEQ\label{eq:min-maxeig-stoch}
\min_{y\in Q} f(y)\equiv\Expect\left[\left\|\pi^{(s)}\left(\sum_{j=1}^p y_j A_j +C\right)\right\|_2\right]-b^Ty
\EEQ
in the variable $y\in\reals^p$ and parameters $A_j\in\symm_n$, for $j=1,\ldots,p$, $b\in\reals^p$ and $C\in \symm_n$, with $1 \leq s \leq n$ controlling the sampling rate, where the function $\|\pi^{(s)}(\sum_{j=1}^p y_j A_j +C)\|_2$ and a subgradient with respect to $y$ are computed using algorithms \ref{alg:matrix-mult} and \ref{alg:low-rank}. For $X\in\symm_n$, we have written $\pi^{(s)}(X)$ the subsampling/scaling operation used in algorithms \ref{alg:matrix-mult} and \ref{alg:low-rank} with
\BEQ\label{eq:pi-sub}
\pi^{(s)}(X)= S,
\EEQ
where $0< s < n$ controls the sampling rate and $S\in\reals^{n \times s}$ is the random matrix defined in algorithm~\ref{alg:low-rank} whose columns are a scaled sample of the columns of X. We will write $S=\pi_{(s)}(X)$ the matrix obtained by subsampling rows as in algorithm~\ref{alg:matrix-mult}. We also define ${\cal A}\in\reals^{n^2 \times p}$ as the matrix whose columns are given by ${\cal A}^{(j)}=\vect(A_j)$, $j=1,\ldots,p$. 


\subsection{Stochastic approximation algorithm}
We show the following lemma approximating the gradient of the function $\|\pi^{(s)}(\sum_{j=1}^p y_j A_j +C)\|_2$ with respect to $y$ and bounding its quadratic variation.

\begin{lemma}\label{lem:grad-var}
Given $A_j\in\symm_n$ with ${\cal A}\in\reals^{n^2 \times p}$ defined as above, for $j=1,\ldots,p$, $b\in\reals^p$, $C\in \symm_n$ and sampling rates $s_1$ and $s_2$, a (stochastic) subgradient of the function $\|\pi^{(s_1)}(\sum_{j=1}^p y_j A_j +C)\|_2-b^Ty$ with respect to~$y$ is given by the vector $w\in\reals^p$ with
\[
w={\cal A}^T\vect(vv^T)u-b
\]
where $v\in\reals^n$ is a leading singular vector of the subsampled matrix $S=\pi^{(s_1)}(\sum_{j=1}^p y_j A_j +C)$ formed in algorithm \ref{alg:low-rank} and $u\in\{-1,1\}$ is the sign of the associated eigenvalue. Furthermore, the product ${\cal A}^T\vect(vv^T)$ can be approximated using algorithm \ref{alg:matrix-mult} to form an approximate gradient
\[
g=\pi^{(s_2)}({\cal A}^T) ~ \pi_{(s_2)}(\vect(vv^T))u-b,
\]
which satisfies
\BEQ\label{eq:quad-var}
\Expect[g]={\cal A}^T\vect(vv^T)u-b\in\partial f(y) \quad\mbox{and}\quad \Expect[\|g\|_2^2]\leq M_*^2 \equiv 2\frac{\|{\cal A}\|_F^2}{s_2}+2\|b\|_2^2.
\EEQ
\end{lemma}
\begin{proof} Iterated expectations give $\Expect[g]=\Expect[w]\in\partial f(y)$. The sampling probabilities $q_i$ used in approximating the matrix vector product ${\cal A}^T\vect(vv^T)$ following algorithm \ref{alg:matrix-mult} are defined as
\[
q_i=\frac{\|{\cal A}_{(i)}\|_2 |\vect(vv^T)_i|}{\sum_{j=1}^{n^2} \|{\cal A}_{(j)}\|_2|\vect(vv^T)_j|}, \quad i=1,\ldots,n^2.
\]
As in \cite[Lemma 3]{Drin07}, the quadratic variation of the approximate product $\pi^{(s_2)}({\cal A}^T) ~ \pi_{(s_2)}(\vect(vv^T))$ is then given by
\[
\Expect[\|\pi^{(s_2)}({\cal A}^T) ~ \pi_{(s_2)}(\vect(vv^T))\|_F^2]=\sum_{i=1}^{n^2} \frac{\|{\cal A}_{(i)}\|_2^2\vect(vv^T)_i^2}{s_2 q_i}.
\]
With $q_i$ defined as above, we get
\BEAS
\sum_{i=1}^{n^2} \frac{\|{\cal A}_{(i)}\|_2^2\vect(vv^T)_i^2}{s_2 q_i} & \leq &  \frac{\left(\sum_{i=1}^{n^2}\|{\cal A}_{(i)}\|_2|\vect(vv^T)_i|\right)^2}{s_2}\\
& \leq & \frac{\|{\cal A}\|_F^2\|vv^T\|_F^2}{s_2}\\
\EEAS
by the Cauchy-Schwarz inequality, because $\|\vect(vv^T)\|^2_2=\|vv^T\|_F^2=\|v\|_2^4=1$, hence the desired result.
\end{proof}

Note that this procedure is not advantageous when $\|{\cal A}^T\|_F \gg \|{\cal A}^T\|_2$ and $s_2$ is small. We now use this result to produce an explicit bound on the complexity of solving problems (\ref{eq:min-maxeig-stoch}) and (\ref{eq:min-maxeig}) by subsampling using a stochastic approximation algorithm. In this section, we let $\|\cdot\|$ be a general norm on $\reals^p$, we write $\|\cdot\|_*$ its dual norm and define $\delta_*(p)$ as the smallest number such that $\|y\|_2\leq \delta_*(p) \|y\|_*$ for all $y\in \reals^p$. Following the notation in \cite[\S 2.3]{Judi07}, we let $\omega(x)$ be a distance generating function, i.e. a function such that
\[
Q^o=\left\{x\in Q:~ \exists y\in \reals^p,~x\in \argmin_{u\in Q} [y^Tu + \omega(u)]\right\}
\]
is a convex set. We assume that $\omega(x)$ is strongly convex on $Q^o$ with modulus $\alpha$ with respect to the norm $\|\cdot\|$, which means
\[
(y-x)^T(\nabla\omega(y)-\nabla\omega(x)) \geq \alpha \|y-x\|^2, \quad x,y\in Q^o.
\]
We then define a prox-function $V(x,y)$ on $Q^o \times Q$ as follows:
\[
V(x,y)\equiv \omega(y) - [ \omega(x)+\nabla \omega(x)^T(y-x)],
\]
which is nonnegative and strongly convex with modulus $\alpha$ with respect to the norm $\|\cdot\|$. The prox-mapping associated to $V$ is then defined as
\BEQ \label{prox-map}
P_x^{Q,\omega}(y) \equiv \argmin_{z\in Q} \{ y^T(z-x) + V(x,z)\}.
\EEQ
Finally, we define the $\omega$ diameter of the set $Q$ as:
\BEQ \label{eq:diameter}
D_{\omega,Q}\equiv(\max_{z\in Q} \omega(z)-\min_{z\in Q} \omega(z))^{1/2}
\EEQ
and we let $\gamma_l$ for $l=0,\ldots,N$ be a step size strategy. 

\begin{algorithm} [h]
\caption{Spectral norm minimization using subsampling} 
\label{alg:stoch-grad} 
\begin{algorithmic}[1]
\REQUIRE Matrices $A_j\in\symm_n$, for $j=1,\ldots,p$, $b\in\reals^p$ and $C\in \symm_n$, sampling rates $s_1$ and $s_2$.
\STATE Pick initial $y_0 \in Q$ 
\FOR{$l=1$ to $N$} 
\STATE Compute $v\in\reals^n$, the leading singular vector of the matrix $\pi^{(s_1)}(\sum_{j=1}^p y_{l,j} A_j +C)$, subsampled according to algorithm \ref{alg:low-rank} with $k=1$ and $s=s_1$.
\STATE Compute the approximate subgradient $g_l=\pi^{(s_2)}({\cal A}^T) ~ \pi_{(s_2)}(\vect(vv^T))-b$, by subsampling the matrix product using algorithm \ref{alg:matrix-mult} and $s=s_2$.
\STATE Set $y_{l+1}=P_{y_l}^{Q,\omega}(\gamma_l g_l)$.
\STATE Update the running average $\tilde y_N= \sum_{k=0}^{N} \gamma_l y_l/\sum_{k=0}^{N}\gamma_l$.
\ENDFOR 
\ENSURE An approximate solution $\tilde y_N\in\reals^p$ of problem (\ref{eq:min-maxeig-stoch}) with high probability.
\end{algorithmic} 
\end{algorithm} 

The following results control the convergence of the robust stochastic approximation algorithm~\ref{alg:stoch-grad} (see  \cite{Judi07}, \cite{Nest09}, \cite{Poly92} or \cite{Nemi83} for further details). We call $\bar y$ the optimal solution of problem~(\ref{eq:min-maxeig-stoch}), the lemma below characterizes convergence speed in expectation.

\begin{lemma}\label{lem:conv-expect}
Given $N>0$, let $M_*$ be defined as in (\ref{eq:quad-var}) by
\[
M_*^2=2\frac{\|{\cal A}\|_F^2}{s_2}+2\|b\|_2^2,
\]
using a fixed step size strategy with 
\[
\gamma_l=\frac{D_{\omega,Q}}{\delta_*(p)M_*}\sqrt{\frac{2}{\alpha N}}, \quad l=1,\ldots,N
\]
we have, after $N$ iterations of algorithm \ref{alg:stoch-grad}
\[
\Expect[f(\tilde y_N)-f(\bar y)] \leq {D_{\omega,Q}}{\delta_*(p)M_*}\sqrt{\frac{2}{\alpha N}}
\]
and 
\[
f(\tilde y_N)-f(\bar y) \geq \epsilon
\]
with probability less than $\frac{D_{\omega,Q}\delta_*(p)M_*}{\epsilon}\sqrt{\frac{2}{\alpha N}}$.
\end{lemma}
\begin{proof}
By construction $\Expect[\|g\|_*^2]\leq \delta_*^2(p) M_*^2$, the rest follows from \cite[\S2.3]{Judi07} for example. 
\end{proof}

Lemma \ref{lem:conv-expect} means that we need at most 
\[
N=\frac{2D_{\omega,Q}^2\delta_*^2(p)M_*^2}{\alpha\epsilon^2\beta^2}
\]
iterations to get an $\epsilon$ solution to problem (\ref{eq:min-maxeig-stoch}) with confidence at least $1-\beta$. Typically, the prox function $\omega$ and the norm are chosen according to the geometry of $Q$, to minimize $N$. The choice of norm also affects $\delta_*(p)$ and obtaining better bounds on $M_*$ in (\ref{eq:quad-var}) for generic norms would further tighten this complexity estimate. 

We now call $y^*$ the solution to the original (deterministic) spectral norm minimization problem~(\ref{eq:min-maxeig}) and bound the suboptimality of~$\tilde y_N$ in the (true) problem~(\ref{eq:min-maxeig}) with high probability.
\begin{theorem}\label{th:conv}
If the sampling rate $s_1$ is set to 
\BEQ\label{eq:opt-sampling}
s_1=\eta^2 \frac{\left\|\sum_{j=1}^p y^*_{j} A_j +C\right\|_2^2}{\epsilon^2 }\NumRank\textstyle\left(\sum_{j=1}^p y^*_{j} A_j +C\right)^2
\EEQ
then after
\BEQ\label{eq:niters}
N=\frac{2D_{\omega,Q}^2\delta_*^2(p)M_*^2}{\alpha\epsilon^2\beta^2}
\EEQ
iterations of algorithm \ref{alg:stoch-grad}, we have
\[
\left\|\sum_{j=1}^p \tilde y_{N,j} A_j +C\right\|_2-b^T\tilde y_N-\left\|\sum_{j=1}^p y^*_{j} A_j +C\right\|_2+b^T y^* \leq 2\epsilon
\]
with probability at least $1-\beta$.
\end{theorem}
\begin{proof}
Recall that we have written $y^*$ the solution to the original (deterministic) problem (\ref{eq:min-maxeig}), $\bar y$~the solution to the approximate (stochastic) problem (\ref{eq:min-maxeig-stoch}) and $\tilde y_N$ the $N$-th iterate of algorithm \ref{alg:stoch-grad} above. Lemma \ref{lem:conv-expect} on the convergence of $\tilde y_N$ to the solution of the stochastic problem in (\ref{eq:min-maxeig-stoch}) means
\[
f(\tilde y_N)-f(\bar y) \leq \epsilon
\]
with probability at least $1-\beta$. By definition, $\bar y$ minimizes the stochastic problem, so in particular $f(\bar y)\leq f(y^*)$, with $f$ the objective value of the {\em stochastic} problem, so we have in fact
\BEQ\label{eq:ineq-opt}
f(\tilde y_N)-f(y^*) \leq \epsilon.
\EEQ
with probability at least $1-\beta$.
Now, with $s_1$ defined as above, Lemma \ref{lem:col-approx} on the quality of the subsampling approximation to $\|.\|_2$ shows that if the sampling rate is set as in (\ref{eq:opt-sampling}) then
\[
\textstyle\Expect\left[\left|\left\|\sum_{j=1}^p y^*_j A_j +C\right\|_2-\left\|\pi^{(s)}\left(\sum_{j=1}^p y_j^* A_j +C\right)\right\|_2\right|\right]\leq \epsilon
\]
and Jensen's inequality yields
\[
\textstyle\left|\left\|\left(\sum_{j=1}^p y^*_j A_j +C\right)\right\|_2-b^Ty^*-f(y^*)\right|\leq \epsilon.
\]
which bounds the difference between the minimum of the (true) problem in (\ref{eq:min-maxeig}) and the value $f(y^*)$ of its stochastic approximation in (\ref{eq:min-maxeig-stoch}), combining this with inequality (\ref{eq:ineq-opt}) we finally get that
\[
\textstyle f(\tilde y_N) - \left\|\left(\sum_{j=1}^p y^*_j A_j +C\right)\right\|_2+b^Ty^*\leq 2\epsilon.
\]
with probability at least $1-\beta$. Applying Jensen's inequality to $\|\cdot\|_2$, using the fact that the subsampling procedure is unbiased, i.e. $\Expect[\pi^{(s)}(X)]=X$ for any $X\in\symm_n$, shows that
\BEAS
\left\|\sum_{j=1}^p \tilde y_{N,j} A_j +C\right\|_2-b^T\tilde y_N & = & \left\|\Expect\left[\pi^{(s)}\left(\sum_{j=1}^p \tilde y_{N,j} A_j +C\right)\right]\right\|_2-b^T\tilde y_{N,j}\\
& \leq & \Expect\left[\left\|\pi^{(s)}\left(\sum_{j=1}^p \tilde y_{N,j} A_j +C\right)\right\|_2\right]-b^T\tilde y_{N,j}\\
& =& f(\tilde y_N)
\EEAS
where the expectation is taken w.r.t. the sampling probability. Hence
\[
\textstyle f(\tilde y_N) - \left\|\left(\sum_{j=1}^p y^*_j A_j +C\right)\right\|_2+b^Ty^*\leq 2\epsilon.
\]
implies
\[
\left\|\sum_{j=1}^p \tilde y_{N,j} A_j +C\right\|_2-b^T\tilde y_N-\left\|\sum_{j=1}^p y^*_{j} A_j +C\right\|_2+b^T y^* \leq 2\epsilon
\]
which is the desired result.
\end{proof}

This result allows us to bound the \emph{oracle} complexity of solving (\ref{eq:min-maxeig}) by subsampling. In practice of course, both the spectral norm and the numerical rank of the solution matrix $\sum_{j=1}^p y^*_{j} A_j +C$ are unknown and guarantees of successful termination (with high probability) depend on the type of {\em stopping criterion} available. Given an exact stopping criterion certifying that $y\in\reals^p$ is optimal (e.g. a target objective value), we can {\em search} for the minimum sampling rate in (\ref{eq:opt-sampling}) by e.g. starting from a low target and doubling the sampling rate until we obtain an optimal solution. On the other hand, if an exact stopping criterion is not available and a more conservative stopping condition is used (e.g. the surrogate duality gap detailed in \S\ref{ss:s-gap}) it is possible for the algorithm to become more expensive than standard subgradient techniques. These two scenarios are detailed below.

\begin{itemize}
\item {\bf Exact stopping criterion.} We assume here that we have a reasonably efficient test for the optimality of the current iterate $\tilde y_{N}$, e.g. a specific target for the objective value below which the optimization procedure can be stopped. In this case, even if we have no a priori knowledge of the rank of the solution matrix, we can search for it, starting with a low guess. The simple lemma below explicitly summarizes the complexity of this procedure.
\begin{lemma}\label{lem:bin-search}
Suppose we start from a sampling rate $s=1$ and run algorithm \ref{alg:stoch-grad} repeatedly, doubling the sampling rate until the stopping criterion certifies the solution is optimal. Then, with probability at least $1-\lceil\log_2 (s_1)\rceil\beta$, algorithm \ref{alg:stoch-grad} needs to be run at most
\[
\lceil\log_2 (s_1)\rceil
\]
times, where $s_1$ is given in (\ref{eq:opt-sampling}), before finding an optimal solution to (\ref{eq:min-maxeig}). 
\end{lemma}
\begin{proof}
Starting from $s=1$, we simply need to double the sampling rate at most $\lceil\log_2 (s_1)\rceil$ before it becomes larger than $s_1$. At the sampling rate $s=s_1$, algorithm \ref{alg:stoch-grad} will produce an optimal solution with prob. $1-\beta$.
\end{proof}

In this scenario, Lemma~\ref{lem:conv-expect} shows that the number of iterations required to reach a target precision $\epsilon$ with confidence greater than $1-\beta \lceil\log_2 (s_1)\rceil$ grows as
\BEQ\label{eq:ntotal}
N^{total}=O\left(\frac{\lceil\log_2 (s_1)\rceil D_{\omega,Q}^2\delta^*(p)^2\left({\|{\cal A}\|_F^2}/{s_2}+\|b\|_2^2\right)}{\alpha\epsilon^2\beta^2}\right)
\EEQ
where  $s_1$ is given in (\ref{eq:opt-sampling}), and the cost of each iteration is detailed in \S\ref{ss:complex} below and the overall complexity of the method is summarized in Table~\ref{tab:complex-stoch}. In fact, we will see in \S\ref{ss:complex} that the complexity of each iteration is dominated by a term $O(sn\log n)$, where $s$ is the sampling rate, and because
\[
\sum_{i=1}^{\lceil\log_2 (s_1)\rceil} 2^i \leq 2^{\lceil\log_2 (s_1)\rceil+1} \leq 4 s_1
\]
we then observe that searching for the minimal sampling rate by repeatedly solving (\ref{eq:min-maxeig-stoch}) for increasing sampling rates will be less than four times as expensive as solving the problem in the oracle case. 

\item {\bf Conservative stopping criterion.} Typically, producing a conservative stopping oracle means computing a surrogate duality gap and we will show in~\S\ref{ss:gap} how this can be done efficiently in some examples. In this case however, it is possible for a conservative stopping criterion to repeatedly fail to detect optimality when searching for the sampling rate $s_1$. The complexity of the subsampling algorithm can then become larger than that of the basic subgradient method.
\end{itemize}

To summarize, when an exact stopping criterion is available, the complexity of finding an optimal solution is equivalent to that of the oracle complexity described in Theorem~\ref{th:conv} and the total number of iterations is bounded by~(\ref{eq:ntotal}). However, when only a conservative stopping condition is available, the algorithm can become more expensive than the classical subgradient method. Note that this early stopping issue is shared by many first-order algorithms, as the theoretical upper bounds available for most first-order methods are usually overly conservative, often by one or two orders of magnitude \cite[\S6]{Nest07}. In these cases too, a conservative stopping criterion is often used to stop the algorithm early. The next section provides a detailed analysis of the complexity of an iteration of algorithm \ref{alg:stoch-grad} as a function of $\epsilon,s_1$ and $s_2$ and the problem data.

\subsection{Complexity}\label{ss:complex}
We now study in detail the complexity of algorithm \ref{alg:stoch-grad}. Suppose we are given a precision target~$\epsilon$ and  fix the sampling rate $s_2$ arbitrarily between 1 and $n^2$, with the sampling rate $s_1$ set as in Theorem~\ref{th:conv}. The cost of each iteration in algorithm \ref{alg:stoch-grad} breaks down as follows.
\begin{itemize}
\item On line 3: Computing the leading singular vector $v$, using algorithm \ref{alg:low-rank} with $k=1$. This means first forming the matrix $(\sum_{j=1}^p y_{l,j} A_j +C)$ and computing the probabilities $q_i$ at a cost of $O(n^2)$ operations. Forming the matrix $S=\pi^{(s_1)}(\sum_{j=1}^p y_{l,j} A_j +C)$ costs $O(ns_1)$ operations. It remains to compute the leading singular vector of $S$ using the Lanczos method at a cost of $O(s_1n\log n)$ (cf. \S\ref{ss:lead-eig} for details). The total numerical cost of this step is then bounded by $c_1 n^2 + c_2 ns_1$ where $c_1$ and $c_2$ are absolute constants. Here, $c_1$ is always less than ten while $c_2$ is the number of iterations required by the Lanczos method to reach a fixed precision target (typically 1e-8 or better here) hence we have $c_1 \ll c_2$.
\item On line 4: Computing the approximate subgradient 
\[
g_l=\pi^{(s_2)}({\cal A}^T) ~ \pi_{(s_2)}(\vect(vv^T))-b,
\]
by subsampling the matrix product using algorithm \ref{alg:matrix-mult}. This means again forming the vector $q$ at a cost of $O(n^2)$ (the row norms of ${\cal A}$ can be precomputed). Computing the subsampled matrix vector product then costs $O(ps_2)$. Both of these complexity bounds have low constants.
\item On line 5: Computing the projection $y_{l+1}=P_{y_l}^{Q,\omega}(\gamma_l g_l)$, whose numerical cost will be denoted by $c(p)$.
\end{itemize}
Let us remark in particular that all $O(n^2)$ operations above only require one pass over the data, which means that the entire data set does not need to fit in memory. Using the bound on the quadratic variation of the gradient computed in Lemma \ref{lem:grad-var}, we can then bound the number of iterations required by algorithm \ref{alg:stoch-grad} to produce a $\epsilon$-solution to problem~(\ref{eq:min-maxeig}) with probability at least $1-\beta$. Let us call $Y^*=\sum_{j=1}^p y^*_{j} A_j +C$, and recall that $\eta=1+\sqrt{8\log(1/\beta)}$, Table \ref{tab:complex-stoch} summarizes these complexity bounds and compares them with complexity bounds for a stochastic approximation algorithm without subsampling.

\begin{table}[H]
\begin{center}
\extrarowheight 1ex
\begin{tabular}{r|c|c}
{\bf Complexity} & Stoch. Approx. & Stoch. Approx. with Subsampling \\
\hline
Per Iter.  & $c_4n^2p+c(p)$ & $ c_2n\log n~\eta^2 \frac{\|Y^*\|_2^2}{\epsilon^2 }\NumRank(Y^*)^2 $ \\
  & & $+c_1n^2 + c_3 p s_2+ c(p)$\\
 &  { } & { } \\
Num. Iter. & $\frac{2D_{\omega,Q}^2\delta^*(p)^2({\|{\cal A}^T\|_2^2}+\|b\|_2^2)}{\alpha\epsilon^2\beta^2}$ & $\frac{2D_{\omega,Q}^2\delta^*(p)^2\left({\|{\cal A}\|_F^2}/{s_2}+\|b\|_2^2\right)}{\alpha\epsilon^2\beta^2}$\\
\end{tabular}
\caption{Complexity of solving problem (\ref{eq:min-maxeig}) using subsampled stochastic approximation method versus original algorithm. Here $c_1,\ldots,c_4$ are absolute constants with $c_1,c_3 \ll c_2,c_4$.\label{tab:complex-stoch}}
\end{center}
\end{table}

We observe that subsampling affects the complexity of solving problem (\ref{eq:min-maxeig}) in two ways. Decreasing the (matrix product) subsampling rate $s_2\in[1,n^2]$ decreases the cost of each iterations but increases the number of iterations in the same proportion, hence has no explicit effect on the total complexity bound. In practice of course, because of higher cache memory speed and better bandwidth on smaller problems, cheaper iterations tend to run more efficiently than more complex ones. Note that the second subsampling step detailed in Lemma~\ref{lem:grad-var} is not advantageous when $\|{\cal A}^T\|_F \gg \|{\cal A}^T\|_2$. When this second subsampling step is skipped the term $\left({\|{\cal A}\|_F^2}/{s_2}+\|b\|_2^2\right)$ is replaced by $\left({\|{\cal A^T}\|_2^2}+\|b\|_2^2\right)$ in the bound on the number of iterations, and the term $c_3ps_2$ becomes $c_3pn^2$ (or less if the matrices are structured).

The impact of the (singular vector) subsampling rate $s_1\in[1,n]$ is much more important however, since computing the leading eigenvector of the current iterate is the most complex step in the algorithm when solving problem (\ref{eq:min-maxeig}) using stochastic approximation. Because $c_1,c_3 \ll c_2$, the complexity per iteration of solving large-scale problems essentially follows
\[
n\log n \,\eta^2 \frac{\|Y^*\|_2^2}{\epsilon^2 }\NumRank(Y^*)^2
\]
hence explicitly depends on both the numerical rank of the solution matrix $Y^*=\sum_{j=1}^p y^*_{j} A_j +C$ and on the relative precision target $\epsilon/\|Y^*\|_2$. This means that problems with simpler solutions will be solved more efficiently than problems whose solutions has a high rank.

The choice of norm $\|\cdot\|$ and distance generating function also has a direct impact on complexity through $c(p)$ and $\delta_*(p)M_*$. Unfortunately here, subsampling error bounds are only available in the Frobenius and spectral norms hence part of the benefit of choosing optimal norm/distance generating function combinations is sometimes lost in the norm ratio bound $\delta_*(p)$. However, choosing a norm/prox function combination according to the geometry of $Q$ can still improve the complexity bound compared to a purely Euclidean setting.

Finally, subsampling could have a more subtle effect on complexity. By construction, solutions to problem (\ref{eq:min-maxeig}) tend to have multiple leading singular values which coalesce near the optimum. Introducing noise by subsampling can potentially break this degeneracy and increase the gap between leading eigenvalues. Since the complexity of the algorithm depends in great part on the complexity of computing a leading singular vector using iterative methods such as the power method or the Lanczos method (cf. Appendix), and the complexity of these methods decreases as the gap between the two leading singular values increases, subsampling can also improve the efficiency of iterative singular value computations. However, outside of simple perturbative regimes, not much is understood at this point about the effect of subsampling on the spectral gap.

\subsection{Surrogate Duality Gap} \label{ss:s-gap}
\label{ss:gap} In practice, we often have no a priori knowledge of $\NumRank(Y^*)^2$ and if the sampling rate~$s$ is set too low, it's possible for the algorithm to terminate at a suboptimal point $Y$ where the subsampling error is less than $\epsilon$ (if the error at the true optimal point $Y^*$ is much larger than $\epsilon$). In order to search for the optimal sampling rate $s$ as in Lemma \ref{lem:bin-search}, we first need to check for optimality in~(\ref{eq:min-maxeig}) and we now show how to track convergence in algorithm \ref{alg:stoch-grad} by computing a surrogate duality gap, at a cost roughly equivalent to that of computing a subgradient. The dual of problem (\ref{eq:min-maxeig}) is written
\BEQ\label{eq:gen-dual}
\BA{ll}
\mbox{maximize} & \Tr(CX) - S_Q(w)\\
\mbox{subject to} & w_j=b_j-\Tr(A_jX),\quad j=1,\ldots,p\\
& \|X\|_{\mathrm{tr}} \leq 1,\\
\EA\EEQ
in the variables $X\in\symm_n$ and $w\in\reals^p$, where $S_Q(v)$ is the support function of the set $Q$, defined as
\[
S_Q(w)\equiv \max_{y\in Q} w^Ty.
\]
For instance, when $Q$ is an Euclidean ball of radius $B$, problem (\ref{eq:gen-dual}) becomes
\BEQ\label{eq:euc-dual}
\BA{ll}
\mbox{maximize} & \Tr(CX) - B \|w\|_2 \\
\mbox{subject to} & w_j=b_j-\Tr(A_jX),\quad j=1,\ldots,p\\
& \|X\|_{\mathrm{tr}} \leq 1,\\
\EA\EEQ
in the variables $X\in\symm_n$ and $w\in\reals^p$. The leading singular vector $v$ in algorithm \ref{alg:stoch-grad} always satisfies $\|vv^T\|_{\mathrm{tr}} \leq 1$, hence we can track convergence in solving (\ref{eq:min-maxeig}) by computing the following surrogate duality gap
\BEQ
\left\|\sum_{j=1}^p y_j A_j +C\right\|_2-b^Ty - v^TCv + S_Q(w)
\EEQ
where $w_j=b_j-v^TA_jv$ for $j=1,\ldots,p$.

\subsection{Minimizing the sum of the $k$ largest singular values} \label{ss:k-sing}
Motivated by applications in statistical learning, we now discuss direct extensions of the results above to the problem of minimizing the sum of the $k$ largest singular values of an affine combination of matrices, written
\BEQ\label{eq:min-ksigv}
\min_{y\in Q}~ \sum_{i=1}^k \textstyle \sigma_i\left(\sum_{j=1}^p y_j A_j +C\right)-b^Ty\\
\EEQ
in the variable $y\in\reals^p$, with parameters $A_j\in\symm_n$, for $j=1,\ldots,p$, $b\in\reals^p$ and $C\in \symm_n$. As in the previous section, we also form its stochastic approximation
\BEQ\label{eq:min-ksigv-stoch}
\min_{y\in Q} f(y)\equiv\Expect\left[\sum_{i=1}^k \sigma_i\left(\pi^{(s)}\left(\sum_{j=1}^p y_j A_j +C\right)\right)\right]-b^Ty
\EEQ
in the variable $y\in\reals^p$, with $1 \leq s \leq n$ controlling the sampling rate. We now prove an analog of Lemma \ref{lem:col-approx} for this new objective function.

\begin{lemma}\label{lem:col-k-approx}
Let $X\in\reals^{m \times n}$ and $\beta\in[0,1]$. Given a precision target $\epsilon>0$, $k\geq 1$ and a matrix $S\in\reals^{m \times s}$ constructed by subsampling the columns of $X$ as in algorithm \ref{alg:low-rank}, let $\eta=1+\sqrt{8\log(1/\beta)}$ and
\BEQ\label{eq:col-samp-k-rate}
s=\eta^2 \frac{(\sum_{i=1}^k \sigma_i(X))^2}{\epsilon^2}\frac{\NumRank(X)^2}{k^2}\kappa(X)^4\Rank(X)
\EEQ
where $\kappa(X)=\sigma_1(X)/\sigma_r(X)$ with $r=\min\left\{k,\Rank(X)\right\}$, we have
\[
\Expect\left[\sum_{i=1}^k \left|\sigma_i(X)-\sigma_i(S)\right|\right]\leq \epsilon
\]
and
\[
\sum_{i=1}^k \left|\sigma_i(X)-\sigma_i(S)\right|\leq \epsilon
\]
with probability at least $1-\beta$.
\end{lemma}
\begin{proof} Because $\Rank(SS^T)\leq \Rank(XX^T)$ by construction, we always have 
\BEAS
\sum_{i=1}^k \left|\sigma_i^2(X)-\sigma_i^2(S)\right| & = & \sum_{i=1}^k \left|\sigma_i(X)-\sigma_i(S)\right|(\sigma_i(X)+\sigma_i(S))\\
& \geq & \sigma_r(X) \sum_{i=1}^k \left|\sigma_i(X)-\sigma_i(S)\right|
\EEAS
where $r=\min\left\{k,\Rank(X)\right\}$. Because the  sum of the $k$ largest singular values is a unitarily invariant norm on $\symm_n$ (see \cite[\S 3.4]{Horn91}), Mirsky's theorem (see \cite[Th. 4.11]{Stew90} for example) shows that
\BEAS
\sum_{i=1}^k \left|\sigma_i^2(X)-\sigma_i^2(S)\right| & = & \sum_{i=1}^k \left|\sigma_i(XX^T)-\sigma_i(SS^T)\right| \\
& \leq & \sum_{i=1}^k \sigma_i(XX^T-SS^T)
\EEAS
and because, by construction, the range of $SS^T$ is included in the range of $XX^T$, we must have $\Rank(XX^T-SS^T)\leq \Rank(XX^T)$ and
\[
\sum_{i=1}^k \sigma_i(XX^T-SS^T) \leq \sqrt{\Rank(X)}~ \|XX^T-SS^T\|_F
\]
Jensen's inequality together with the matrix multiplication result in Lemma \ref{ref:lem-col-sample} yield
\[
\Expect[\|SS^T-XX^T\|_F] \leq \displaystyle \frac{\|X\|_F^2}{\sqrt{s}}
\]
and
\[
\|SS^T-XX^T\|_F \leq \frac{\eta \|X\|_F^2}{\sqrt{s}}
\]
with probability at least $1-\beta$. Combining these inequalities with the sampling rate in~(\ref{eq:col-samp-k-rate}) 
\[
s=\eta^2 \frac{\|X\|_F^4\Rank(X)}{\epsilon^2 \sigma_r(X)^2}
\]
and using
\[
\frac{\|X\|_F^4}{(\sum_{i=1}^k \sigma_i(X))^2 \sigma_r(X)^2} \leq \frac{\NumRank(X)^2}{k^2}\kappa(X)^4
\]
yields the desired result.
\end{proof}

Once again, the subsampling rate in the above lemma has a clear interpretation, 
\[
s_1=\eta^2 \frac{(\sum_{i=1}^k \sigma_i(X))^2}{\epsilon^2}\frac{\NumRank(X)^2}{k^2}\kappa(X)^4\Rank(X)
\]
is the product of a term representing relative precision, a term reflecting the rank of $X$ and a term in $\kappa(X)$ representing its (pseudo) condition number. Note that the bound can be further refined when $\sigma_r \leq \epsilon$. Lemma \ref{lem:col-k-approx} allows us to compute the gradient by subsampling when using algorithm \ref{alg:stoch-grad} to solve problem (\ref{eq:min-ksigv}). The remaining steps in the algorithm are identical, except that the matrix $vv^T$ is replaced by a combination of matrices formed using the $k$ leading singular vectors, with Frobenius norm $\sqrt{k}$. The cost of each iteration is dominated by the term 
\[
c_2 k s_1 m\log n \, + c_1 nm + c_3 p s_2 + c(p)
\]
with $s_1$ defined above, and the total number of iterations growing as
\[
\frac{2D_{\omega,Q}^2\delta^*(p)^2\left({k \|{\cal A}\|_F^2}/{s_2}+\|b\|_2^2\right)}{\alpha\epsilon^2\beta^2}.
\]
When $\|{\cal A}^T\|_F \gg \|{\cal A}^T\|_2$ and $s_2$ is small, the second subsampling step detailed in Lemma~\ref{lem:grad-var} is skipped and the cost per iteration becomes
\[
c_2 k s_1 m\log n \, + c_1 nm + c_3 p mn + c(p)
\]
and the total number of iterations grows as
\[
\frac{2D_{\omega,Q}^2\delta^*(p)^2\left({k \|{\cal A}^T\|_2^2}+\|b\|_2^2\right)}{\alpha\epsilon^2\beta^2}
\]
because $M_*^2 \leq 2\left({k\|{\cal A}^T\|_2^2}+\|b\|_2^2\right)$ in this case.

\section{Applications \& numerical results}
\label{s:numexp} In this section, we first detail a few instances of problem (\ref{eq:min-maxeig}) arising in statistical learning. We then study the numerical performance of the methods detailed here on large scale problems.

\subsection{Spectral norm minimization}
\label{ss:nrom-min} For a given matrix $A\in\symm_n$, we begin by studying a simple instance of problem (\ref{eq:min-maxeig}) written
\BEQ\label{eq:min-spca}
\BA{ll}
\mbox{minimize} & \|A+U\|_2\\
\mbox{subject to} & |U_{ij}| \leq \rho,\quad i,j=1,\ldots,n\\
\EA\EEQ
in the matrix $U\in\symm_n$. This problem is closely related to a relaxation for sparse PCA (see \cite{dAsp04a}) and we use it in the next section to test the numerical performance of algorithm~\ref{alg:stoch-grad}. The complexity of the main step in the algorithm (i.e. computing the gradient) is controlled by the sampling rate in Lemma~\ref{lem:col-approx}, which is written
\[
s_1=\eta^2 \frac{\|A+U^*\|_2^2}{\epsilon^2 }\NumRank(A+U^*)^2
\]
where $U^*\in\symm_n$ is the optimal solution to problem (\ref{eq:min-spca}). The prox function used here is the square Euclidean norm, and the prox-mapping is then a simple Euclidean projection on the box $[-\rho,\rho]^{n^2}$. The cost of each iteration is then dominated by the term 
\[
c_2 s_1 n\log n \, + (c_1+c_3+1) n^2.
\]
with $s_1$ defined above, and the total number of iterations grows as
\[
\frac{4\lceil\log_2 (s_1)\rceil n^2\rho^2}{\epsilon^2\beta^2},
\]
because the gradient always has norm one in this problem and the second subsampling step in Lemma~\ref{lem:grad-var} is not beneficial.

\subsection{Matrix factorization and collaborative filtering}
\label{ss:matrix-fact} Matrix factorization methods have been heavily used to solve collaborative filtering problems (e.g. the {\em Netflix} problem) and we refer the reader to \cite{Sreb04}, \cite{Bach07}, \cite{Rech07} or \cite{Cand08} for details. \cite{Sreb04} focuses on the following problem instance 
\BEQ\label{eq:min-rank-hinge}
\mbox{minimize}\quad \|X\|_\mathrm{tr} + c \sum_{(i,j)\in S}\max(0,1-X_{ij}M_{ij})
\EEQ
in the variable $X\in \reals^{m \times n}$, where $M$ is a sparse matrix of ratings, $S$ is the set of known ratings (typically small), and $c>0$ is a parameter controlling the rank versus accuracy tradeoff. Here, the trace norm can be understood as a convex lower bound on the rank function (as in \cite{Boyd00}) but sometimes also has a direct interpretation in terms of learning (see \cite{Sreb04}). The dual of this problem is written 
\[\BA{ll}
\mbox{maximize} & \sum_{ij} Y_{ij}\\
\mbox{subject to} & \|Y \circ M\|_2 \leq 1\\
& 0 \leq Y_{ij} \leq c
\EA\]
in the variable $Y\in\reals^{n \times m}$, where $Y \circ M$ is the Schur (componentwise) product of $Y$ and $M$. Because $M$ is usually sparse, this problem is typically sparse too (i.e. most of the coefficients of $Y$ can be set to $c$). This last problem is equivalent to
\[\BA{ll}
\mbox{minimize} & \|Y \circ M\|_2 \\
\mbox{subject to} & \sum_{ij} Y_{ij}= 1\\
& 0 \leq Y_{ij} \leq d
\EA\]
for some $d>0$, which is a particular instance of~(\ref{eq:min-maxeig}). In the numerical experiments that follow, we focus on a simpler formulation of (\ref{eq:min-ksigv}) written
\BEQ\BA{ll}\label{eq:min-tracenorm}
\mbox{minimize} & \sum_{i=1}^k \textstyle \sigma_i (X)\\
\mbox{subject to} & X_{ij}=M_{ij},\quad (i,j)\in S\\
& |X_{ij}| \leq B
\EA\EEQ
in the variable $X \in \reals^{m \times n}$, for some $k\in[1,n]$ and $B>0$, which is also a particular instance of problem (\ref{eq:min-ksigv}). We assume $n\geq m$. This replaces the hinge-loss penalty in problem~(\ref{eq:min-rank-hinge}) with equality constraints. The classic trace heuristic uses $\|X\|_\mathrm{tr}$ instead of $\sum_{i=1}^k  \sigma_i (X)$ but we will see in what follows that minimizing the later term also tends to produce low rank solutions.

In this particular case, the complexity of the main step in the algorithm (i.e. computing the gradient) is controlled by the sampling rate in Lemma~\ref{lem:col-k-approx}, which can be simplified here to
\[
s_1=\eta^2 \frac{ \left\|Y^*\right\|_\mathrm{tr}^2}{\epsilon^2}\kappa(Y^*)^2\Rank(Y^*)
\]
where $Y^*=\sum_{j=1}^p y^*_j A_j +C$ and $\kappa(Y^*)=\sigma_1(Y^*)/\sigma_r(Y^*)$ with $r=\Rank(Y^*)$. The prox function used here is the square Euclidean norm, and the prox-mapping is then a simple Euclidean projection on the box $[-B,B]^{m \times n}$ for the coefficients whose rating is not given. As in~(\ref{ss:k-sing}), the cost of each iteration is  dominated by the term 
\[
c_2 s_1 m\log n \, + (c_1+c_3+1) nm.
\]
with $s_1$ defined above, and the total number of iterations growing as
\[
\frac{4\lceil\log_2 (s_1)\rceil nm B^2 k}{\epsilon^2\beta^2},
\]
because the gradient has Frobenius norm $\sqrt{k}$ in this case. This bound can be further refined when $\sigma_r \leq \epsilon$. In practice, the complexity of solving problem~(\ref{eq:min-tracenorm}) can often be further reduced using the simple observation that an optimal solution of (\ref{eq:min-ksigv}) will also be optimal in (\ref{eq:min-tracenorm}) whenever $\Rank(Y^*_k) < k$, where $Y^*_k$ is the optimal solution to (\ref{eq:min-ksigv}) here. Once again, the sampling rate $s$ has a natural interpretation as the product of a relative precision term, a term reflecting the condition number of the solution and the rank of the optimal solution. It means in particular that problems whose solutions have a lower rank are explicitly easier to solve than problems with more complex (higher rank) solutions. Of course, much faster specialized algorithms are available for this problem, but those methods that exploit the fact that the solution is low-rank (like alternating minimization) are non-convex.

\subsection{Numerical experiments}
\label{s:numres}
In this section, we test the quality of the subsampling approximations detailed in Section~\ref{s:random} on various matrices. We also evaluate the performance of the algorithms detailed above on large scale problem instances. Numerical code reproducing these experiments is available from the author's webpage.

\paragraph{Randomized low-rank approximations.}
Here, we first measure the quality of the randomized low-rank matrix approximation on both randomly generated matrices and on covariance matrices formed using gene expression data. Because the spectrum of naive large scale random matrices is very structured, these examples are too simple to appropriately benchmark numerical error in algorithm \ref{alg:low-rank}. Fortunately, as we will see below, generating random symmetric matrices with a given spectral measure is straightforward.

Suppose $X\in\symm_n$ is a matrix with normally distributed coefficients, $X_{ij}\sim\mathcal{N}(0,1)$, $i,j=1,\ldots,n$. If we write its QR decomposition, $X=QR$ with $Q,~R\in \reals^{n \times n}$, then the orthogonal matrix $Q$ is Haar distributed on the orthogonal group $\mathcal{O}_n$ (see \cite{Diac03} for example). This means that to generate a random matrix with given spectrum $\mu\in\reals^n$, we generate a normally distributed matrix $X$, compute its QR decomposition and the matrix $Q\diag(\mu)Q^T$ will be uniformly distributed on the set of symmetric matrices with spectrum $\mu$. Because the spectral measure of ``natural'' covariance matrices often follows a power law (Tracy-Widom in the Gaussian case, see \cite{John01} and \cite{El-K07} for a discussion), we sample the spectrum $\mu$ from a beta distribution with various exponents to get realistic random matrices with a broad range of numerical ranks. We also use a covariance matrix formed using the gene expression data set in \cite{Alon99}.

\begin{figure}[hp]
\begin{center}
\begin{tabular}{cc}
\psfrag{numrank}[t][b]{$\NumRank(X)$}
\psfrag{relerr}[b][t]{$\epsilon/\|X\|_2$}
\includegraphics[width=0.49 \textwidth]{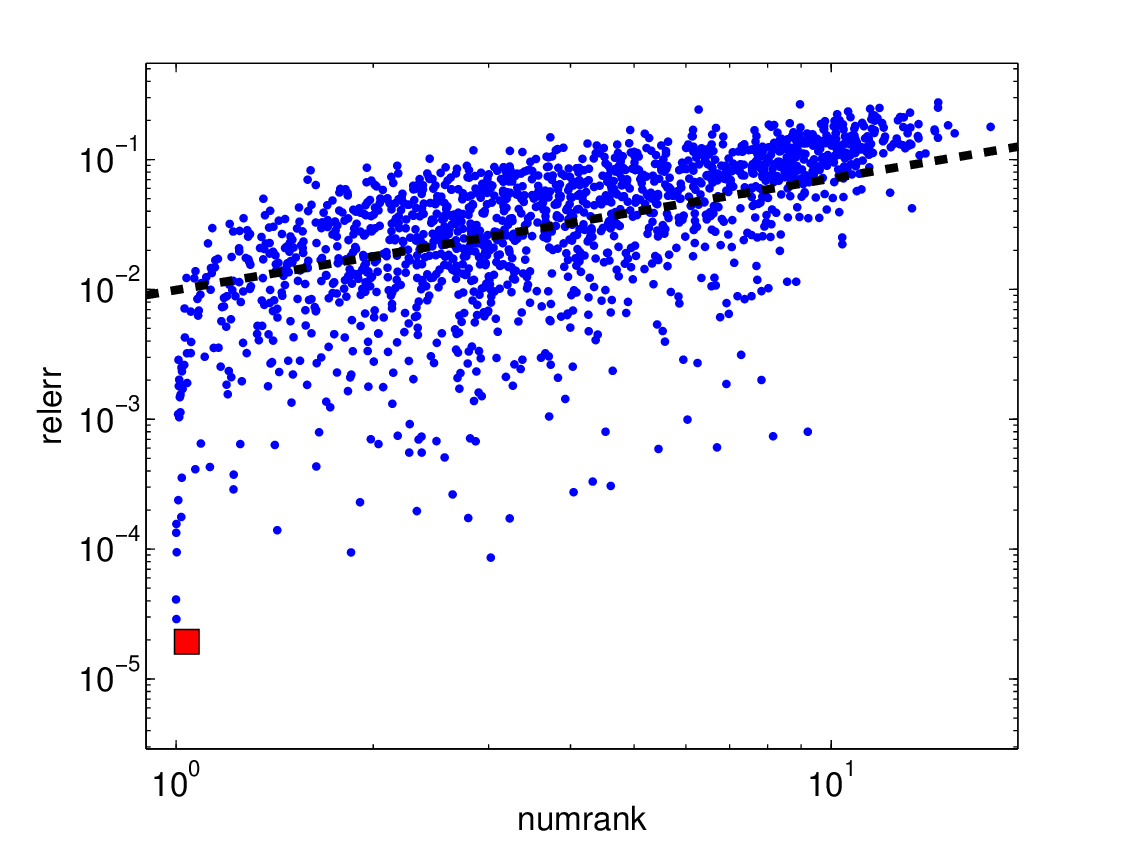}&
\psfrag{erratio}[t][b]{Error / Theoretical error}
\psfrag{occur}[b][t]{\# occurences}
\includegraphics[width=0.45\textwidth]{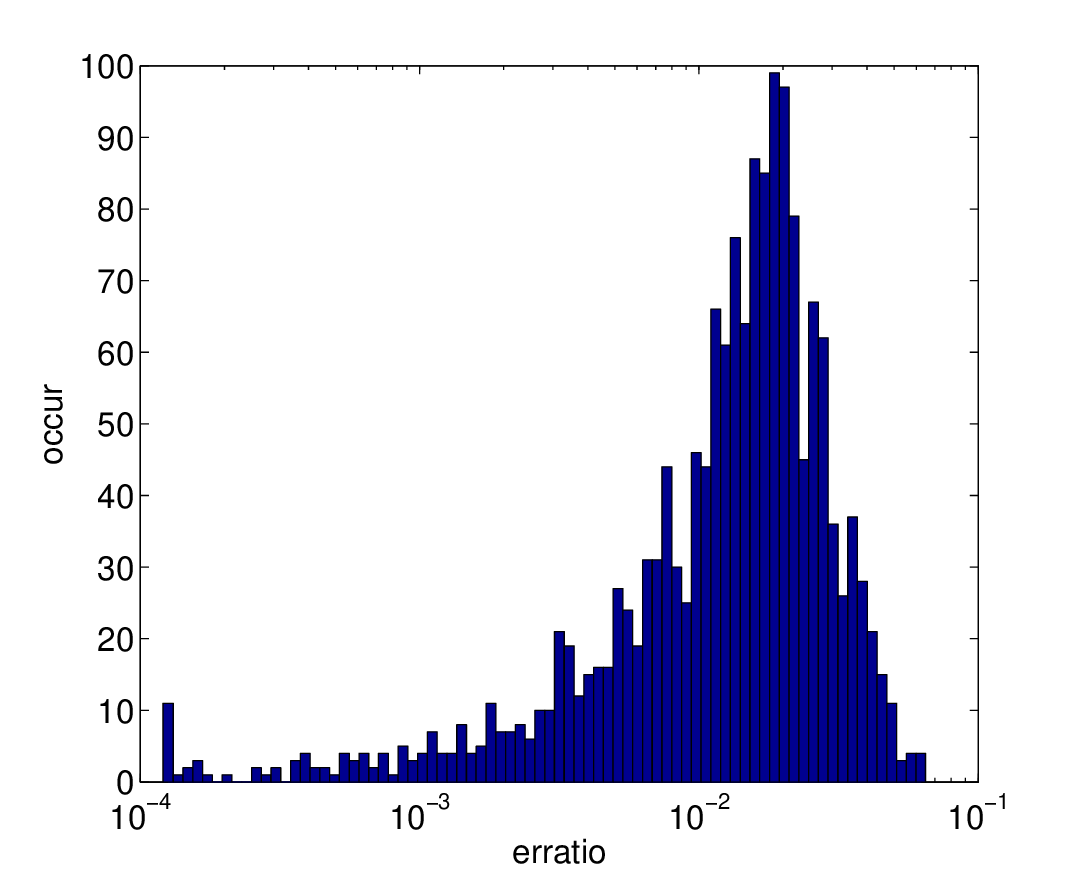}
\end{tabular}
\caption{\textit{Left:} Loglog plot of relative error $\epsilon/\|X\|_2$ versus numerical rank $\NumRank(X)$ with 20\% subsampling and $n=500$ on random matrices (blue dots) and gene expression covariance (red square). The dashed line has slope one in loglog scale. \textit{Right:} Histogram plot in semilog scale of relative error $\epsilon/\|X\|_2$ over theoretical bound $\eta\NumRank(X)/\sqrt{s}$ for random matrices with $n=500$.
\label{fig:err-vs-rank}}
\end{center}
\end{figure}

\begin{figure}[hp]
\begin{center}
\begin{tabular}{cc}
\psfrag{srate}[t][b]{Sampling rate $s$}
\psfrag{sqrelerr}[b][t]{$\epsilon/\|X\|_2$}
\includegraphics[width=0.49 \textwidth]{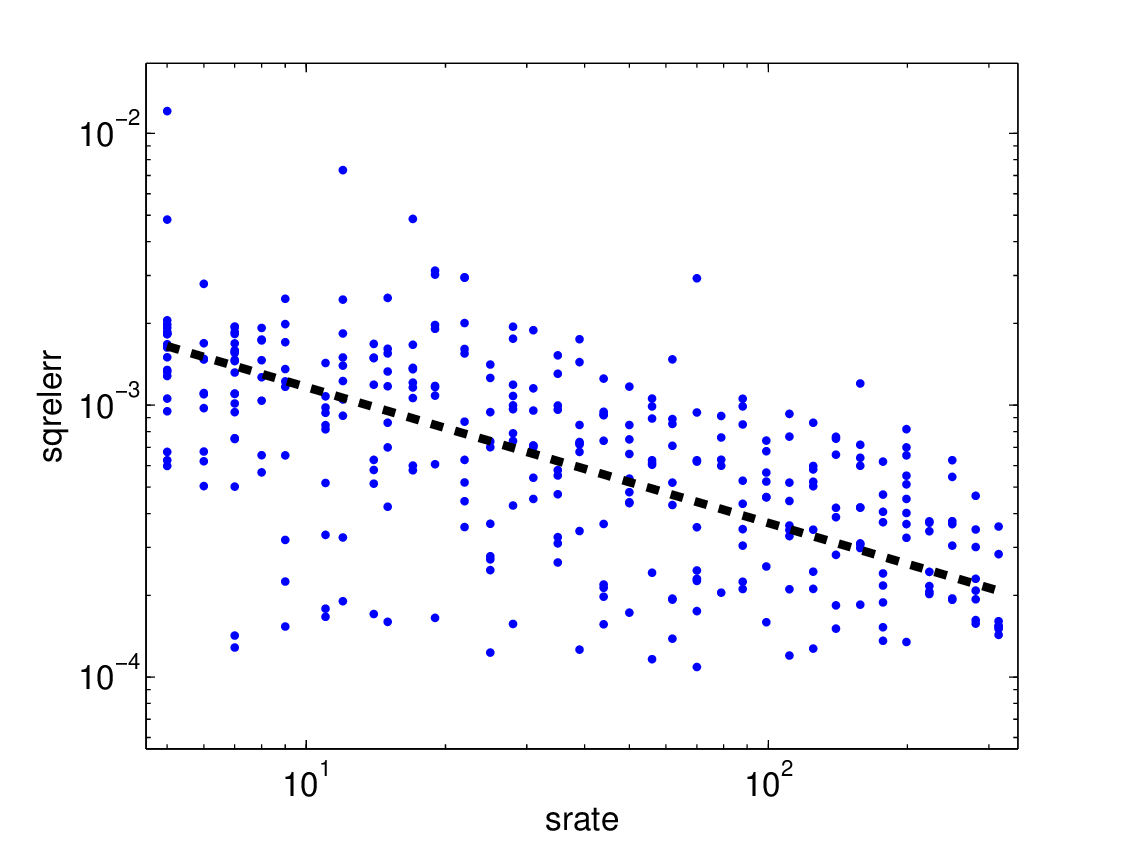}&
\psfrag{sratio}[t][b]{Sampling ratio $s/n$}
\psfrag{speedup}[b][t]{Speedup factor}
\includegraphics[width=0.49\textwidth]{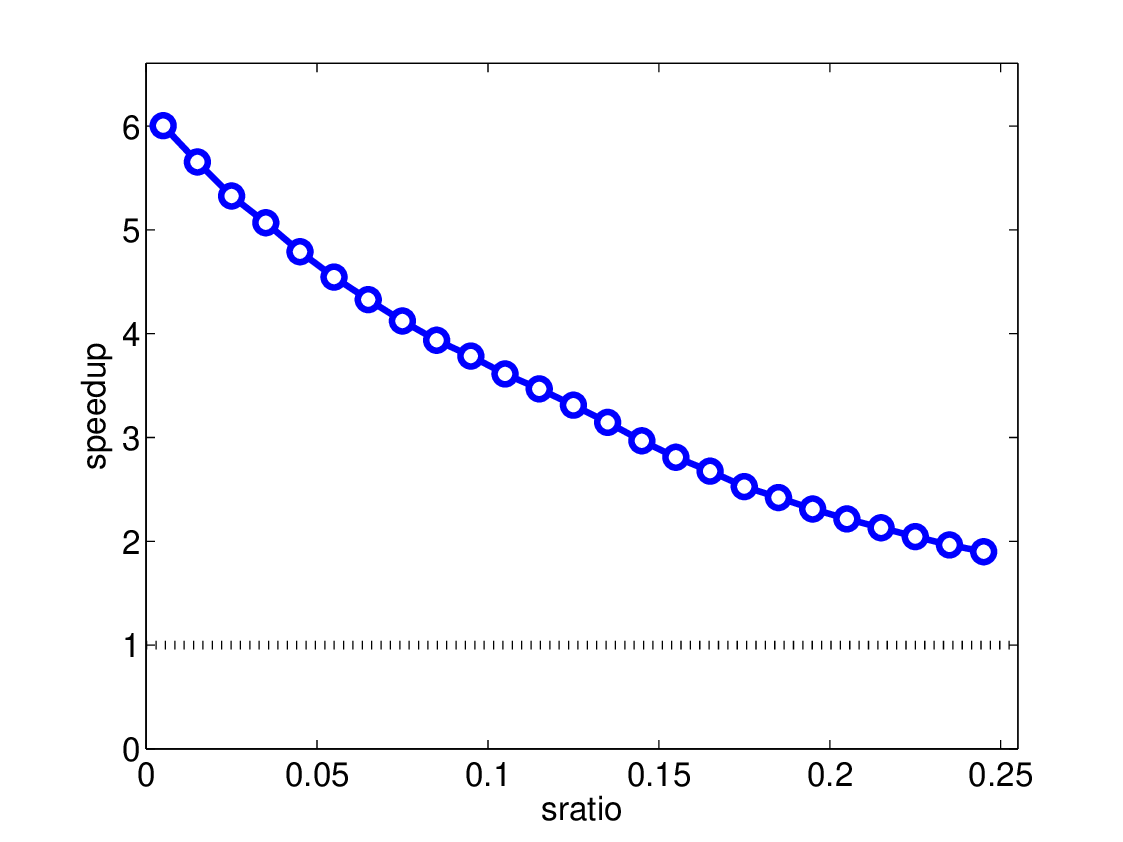}
\end{tabular}
\caption{\textit{Left:} Loglog plot of relative error $\epsilon/\|X\|_2$ versus sampling rate $s$ for a gene expression covariance  with $n=500$. The dashed line has slope -1/2 in loglog scale. \textit{Right:} Plot of median speedup factor in computing largest magnitude eigenvalue, using ARPACK with and without subsampling on a gene expression covariance matrix with $n=2000$, for various values of the sampling ratio $s/n$.
\label{fig:err-vs-sample}}
\end{center}
\end{figure}

In Figure \ref{fig:err-vs-rank}, we plot relative error $\epsilon/\|X\|_2$ against the numerical rank $\NumRank(X)$ in loglog scale with 20\% subsampling and $n=500$ on random matrices generated as above and on the gene expression covariance from \cite{Alon99}. We notice that, on these experiments, the relative error grows at most linearly with the numerical rank of the matrix, as predicted by Lemma~\ref{lem:col-approx}. We then plot the histogram in semilog scale of relative error $\epsilon/\|X\|_2$ over theoretical bound $\eta\NumRank(X)/\sqrt{s}$ for random matrices with $n=500$. In Figure \ref{fig:err-vs-sample}, we plot relative error $\epsilon/\|X\|_2$ versus sampling rate $s$, in loglog scale, for a gene expression covariance with $n=500$. Once gain, the error decreases as $1/\sqrt{s}$ as predicted by Lemma~\ref{lem:col-approx}. We also plot the median speedup factor (over ten runs) in computing largest magnitude eigenvalues using ARPACK with and without subsampling on a gene expression covariance matrix with $n=2000$, for various values of the sampling ratio $s/n$. Note that both exact and subsampled eigenvalues are computed using direct MEX calls to ARPACK by \cite{Leho98}, as \texttt{eigs} (MATLAB's interface to ARPACK) carries a massive overhead. In all the experiments above, the confidence level used in computing $\eta$ was set to 99\%.

\begin{figure}[h]
\begin{center}
\begin{tabular}{cc}
\psfrag{cpu}[t][b]{CPU time (secs.)}
\psfrag{oval}[b][t]{Objective value}
\includegraphics[width=0.49 \textwidth]{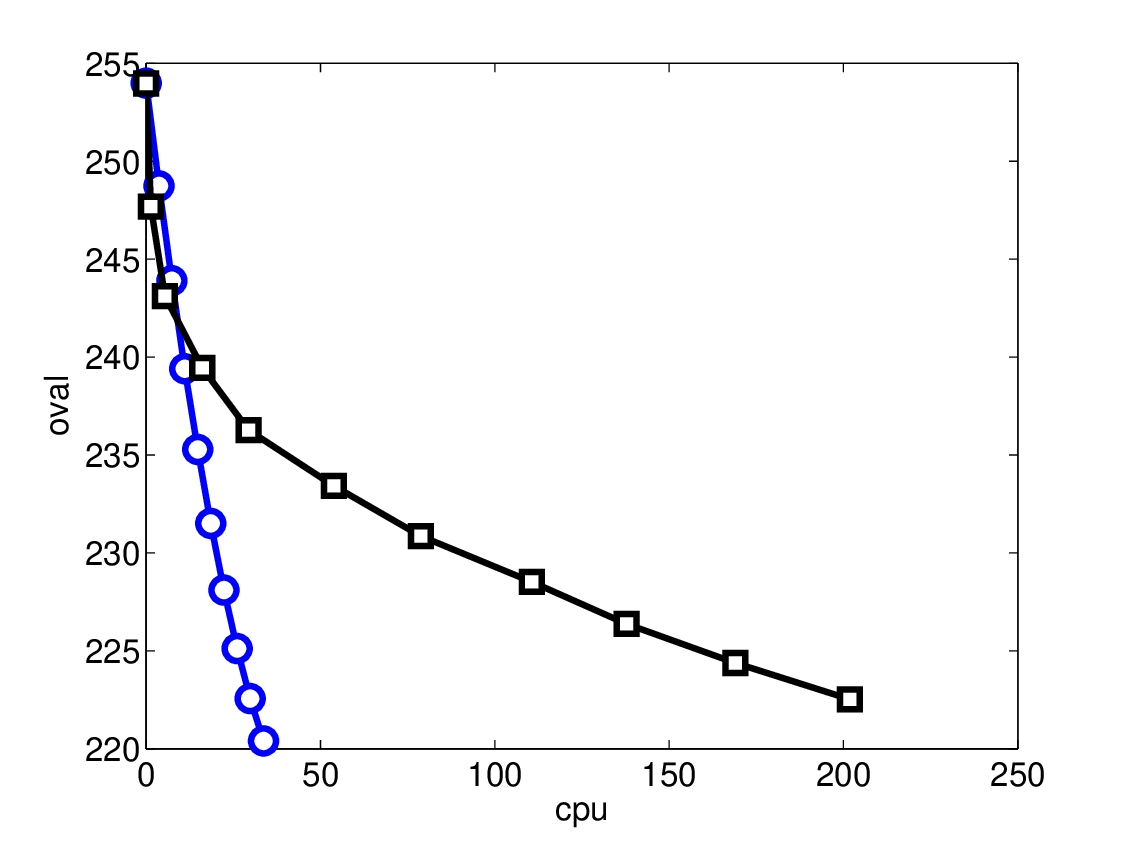}&
\psfrag{cpu}[t][b]{CPU time (secs.)}
\psfrag{gap}[b][t]{Surrogate gap}
\includegraphics[width=0.49\textwidth]{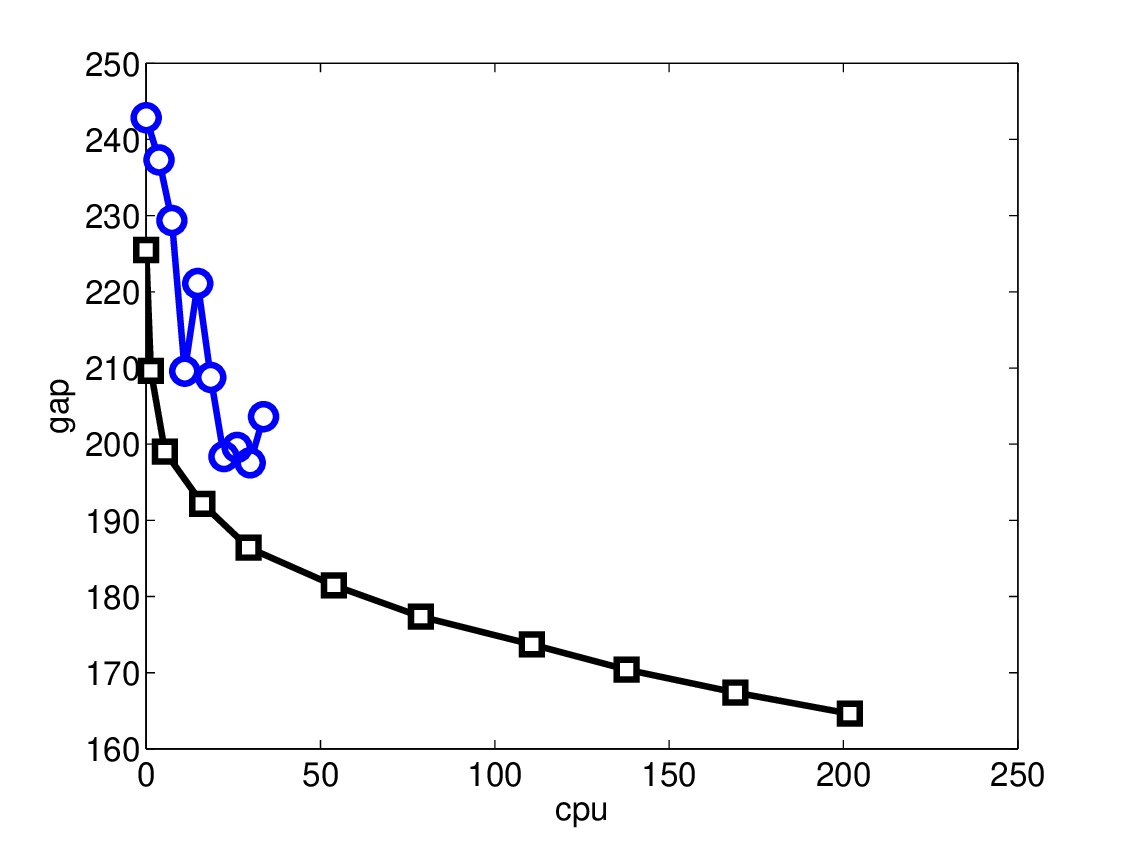}
\end{tabular}
\caption{\textit{Left:} Objective value versus CPU for a sample matrix factorization problem in dimension 100, using a deterministic gradient (squares) or a subsampled gradient with subsampling rate set at 20\% (circles). \textit{Right:} Surrogate duality gap versus CPU time on the same example.
\label{fig:sample-trace}}
\end{center}
\end{figure}

\paragraph{Stochastic approximation with subsampling.} In Figure \ref{fig:sample-trace}, we generate a sample ratings matrix $X=VV^T$ for the collaborative filtering problem~(\ref{eq:min-tracenorm}) in \S\ref{ss:matrix-fact}, where $V$ is a discrete feature matrix $V\in[0,4]^{n \times 3}$, with $n=100$. We ``observe'' only 30\% of the coefficients in $V$ and solve problem (\ref{eq:min-tracenorm}) with $k=4$ and $B=10$ to approximately reconstruct the full ratings matrix. We plot objective value versus CPU time in seconds for this sample matrix factorization problem, using a stochastic approximation algorithm with deterministic gradient or the subsampled gradient algorithm \ref{alg:stoch-grad} with subsampling ratio~$s_1/n$ set at 20\%. We also plot surrogate duality gap versus CPU time on the same example. We notice that while the subsampled algorithm converges much faster than the deterministic one, the quality of the surrogate dual points and duality gap  produced using subsampled gradients as in \S\ref{ss:gap} is worst than in the deterministic case. 

In Table \ref{tab:box-spectral}, using the same 20\% sampling rate we compare CPU time versus problem dimension $n$ for subsampled and deterministic algorithms when solving the following instance of problem (\ref{eq:min-maxeig}) 
\[\BA{ll}
\mbox{minimize} & \|C+X\|_2\\
\mbox{subject to} & \|X\|_\infty \leq \rho
\EA\]
in the variable $X\in\symm_n$ where $C$ is a covariance matrix constructed using the $n$ variables in the gene expression data set of \cite{Alon99} with maximum variance, for various values of $n$, and $\rho=1/2$. We run 200 iterations of the deterministic algorithm and run the stochastic algorithm until it reaches the best value found by the deterministic method. Finally, using random ratings data generated as above, we solve sample collaborative filtering solve problems (\ref{eq:min-tracenorm}) with $k=4$ and $B=10$ for ratings matrix of various dimensions $n$. 
We run 5000 iterations of the deterministic algorithm and run the stochastic algorithm until it reaches the best value found by the deterministic method.

Of course, specialized software packages solve much larger problems and this example is only here to illustrate the potential speedup. We report median CPU time over ten sample problems in Table \ref{tab:collab-filter}. Here, subsampling speeds up the algorithm by an order of magnitude, however the stochastic approximation algorithm is still not competitive with (non convex) local minimization techniques over low rank matrices.

\begin{table}[H]
\begin{center}
\begin{tabular}{r|c|c|c}
$n$ & Deterministic & Subsampling  & Speedup factor\\
\hline
500 & 5 & 5 & 0.92 \\
750 & 19 & 13 & 1.40\\
1000 & 32 & 24 & 1.31 \\
1500 & 107 & 58 & 1.84 \\
2000 & 281 & 120 & 2.34   
\end{tabular}
\caption{CPU time (in seconds) versus problem dimension $n$ for deterministic and subsampled stochastic approximation algorithms on spectral norm minimization problems. \label{tab:box-spectral}}
\end{center}
\end{table}

\begin{table}[H]
\begin{center}
\begin{tabular}{r|c|c|c}
$n$ & Deterministic & Subsampling  & Speedup factor\\
\hline
100  & 154 & 23 & 6.67 \\
200  & 766 & 63 & 12.2 \\
500  & 4290 & 338 & 12.7   
\end{tabular}
\caption{Median CPU time (in seconds) versus problem dimension $n$ for deterministic and subsampled stochastic approximation algorithms on collaborative filtering problems. \label{tab:collab-filter}}
\end{center}
\end{table}

\section{Appendix}
The complexity results detailed above heavily rely on the fact that extracting {\em one} leading eigenvector of a symmetric matrix $X\in\symm_n$ can be done by computing a few matrix vector products. While this simple fact is easy to prove using the power method when the eigenvalues of $X$ are well separated, the problem becomes significantly more delicate when the spectrum of $X$ is clustered. The section that follows briefly summarizes how modern numerical methods solve this problem in practice.

\subsection{Computing one leading eigenvector of a symmetric matrix} \label{ss:lead-eig}
We start by recalling how packages such as LAPACK \cite{Ande99} form a full eigenvalue (or Schur) decomposition of a symmetric matrix $X\in\symm_n$. The algorithm is strikingly stable and, despite its $O(n^3)$ complexity, often competitive with more advanced techniques when the matrix $X$ is small. We then discuss the problem of approximating one leading eigenpair of $X$ using Krylov subspace methods with complexity growing as $O(n^2\log n)$ with the dimension (or less when the matrix is structured).

\paragraph{Full eigenvalue decomposition.}
Full eigenvalue decompositions are computed by first reducing the matrix $X$ to symmetric tridiagonal form using Householder transformations, then diagonalizing the tridiagonal factor using iterative techniques such as the QR or divide and conquer methods for example (see \cite[Chap. 3]{Stew01} for an overview). The classical QR algorithm (see \cite[\S8.3]{Golu90}) implicitly relied on power iterations to compute the eigenvalues and eigenvectors of a symmetric tridiagonal matrix with complexity $O(n^3)$, however more recent methods such as the MRRR algorithm by \cite{Dhil03} solve this problem with complexity $O(n^2)$. Starting with the third version of LAPACK, the MRRR method is part of the default routine for diagonalizing a symmetric matrix and is implemented in the \texttt{STEGR} driver (see \cite{Dhil06}). 

Overall, the complexity of forming a {\em full} Schur decomposition of a symmetric matrix $X\in\symm_n$ is then $4n^3/3$ flops for the Householder tridiagonalization, followed by $O(n^2)$ flops for the Schur decomposition of the tridiagonal matrix using the MRRR algorithm.

\paragraph{Computing one leading eigenpair.} We now give a brief overview of the complexity of computing leading eigenpairs using Krylov subspace methods and we refer the reader to \cite[\S4.3]{Stew01}, \cite[\S8.3, \S9.1.1]{Golu90} or \cite{Saad92} for a more complete discussion. Let $u\in\reals^n$ be a given vector, we form the following {\em Krylov} sequence
\[
\left\{u,Xu,X^2u,\ldots,X^ku\right\}
\]
by computing $k$ matrix vector products. If we call ${\cal K}_k(X,u)$ the subspace generated by these vectors and write $X=\sum_{i=1}^n \lambda_i x_ix_i^T$ a spectral decomposition of $X$, assuming, for now, that 
\[
\lambda_1 > \lambda_2 \geq \ldots \geq \lambda_n,
\]
one can show using Chebyshev polynomials (see e.g. \cite[\S4.3.2]{Stew01} for details) that
\[
\tan \angle\left(x_1,{\cal K}_k(X,u)\right) \lesssim \frac{\tan\angle(x_1,u)}{\left(1+2\sqrt{\eta+\eta^2}\right)^{k-1}}
\quad \mbox{where} \quad \eta=\frac{\lambda_1-\lambda_2}{\lambda_2-\lambda_n},
\]
in other words, after a few iterations, Krylov subspaces contain excellent approximations of leading eigenpairs of $X$. 

This result is exploited by the Lanczos procedure to extract approximate eigenpairs of $X$ called {\em Ritz} pairs (see \cite[Chap. 9]{Golu90} or \cite[\S5.1.2]{Stew01} for a complete discussion). In practice, the matrix formed by the Krylov sequence is very ill-conditioned (as $X^ku$ gets increasingly close to the leading eigenvector) so the Lanczos algorithm simultaneously updates an orthogonormal basis for ${\cal K}_k(X,u)$ and a {\em partial} tridiagonalization of $X$. The Lanczos procedure is described in Algorithm~\ref{alg:lanczos} and requires $k$ matrix vector products and an additional $4nk$ flops. Note that the {\em only} way in which the data in $X$ is accessed is through the matrix vector products $Xu_j$.

\begin{algorithm}[h!]
\caption{Lanczos decomposition.} 
\label{alg:lanczos} 
\begin{algorithmic}[1]
\REQUIRE Matrices $X\in\symm_n$ and initial vector $u_1\in\reals^n$.
\STATE Set $u_0=0$ and $\beta_0=0$.
\FOR{$j=1$ to $k$} 
\STATE Compute $v=Xu_j$.
\STATE Set $\alpha_j=u_j^Tv$.
\STATE Update $v=v-\alpha_ju_j-\beta_{j-1}u_{j-1}$.
\STATE Set $\beta_j=\|v\|_2$.
\STATE Set $u_{j+1}=v/\beta_j$.
\ENDFOR 
\ENSURE A Lanczos decomposition
\[
XU_k=U_kT_k+\beta_ku_{k+1}e^T_{k},
\]
where $U_k\in\reals^{n\times k}$ is orthogonal and $T_k\in\symm_k$ is symmetric tridiagonal.
\end{algorithmic} 
\end{algorithm} 

In theory, one could then diagonalize the matrix $T_k$ (which costs $O(k^2)$ using the MRRR algorithm as we have seen above) to produce Ritz vectors. In practice, key numerical difficulties often arise. First, finite precision arithmetics cause a significant loss of orthogonality in $U_k$. This is remedied by various reorthogonalization strategies (cf. \cite[\S5.3.1]{Stew01}). A more serious problem is clustered or multiple eigenvalues in the spectrum periphery. In fact, it is easy to see that Krylov subspace methods cannot isolate multiple eigenvalues. Assume that the leading eigenvalue has multiplicity two for example, we then have
\[
A^ku=((x_1^Tu)x_1 + (x_2^Tu) x_2) \lambda_1^k + (x_3^Tu) x_3 \lambda_3^k +\ldots + (x_n^Tu) x_n \lambda_n^k
\]
and everything happens as if the eigenvalue $\lambda_1$ was simple and the matrix $X$ had a larger nullspace. This is not a problem in the optimization problems discussed in this paper, since we need only {\em one} eigenvector in the leading invariant subspace, not the entire eigenspace.

Clustered eigenvalues (i.e. a small gap between the leading eigenvalue and the next one, not counting multiplicities) are much more problematic. The convergence of Ritz vectors cannot be established by the classical Chebyshev bounds described above, and various references provide a more refined analysis of this scenario (see \cite{Parl82}, \cite{Van-87}, \cite{Kucz92} among others). Successful termination of a {\em deterministic} Lanczos method can never be guaranteed anyway, since in the extreme case where the starting vector is orthogonal to the leading eigenspace, the Krylov subspace contains no information about leading eigenpairs. In practice, Lanczos solvers use {\em random} initial points. In particular, \cite[Th.4.2]{Kucz92} show that, for any matrix $X\in\symm_n$ (including matrices with clustered spectrum), starting the algorithm at a random $u_1$ picked uniformly over the sphere means the Lanczos decomposition will produce a leading Ritz pair with {\em relative} precision $\epsilon$ in
\[
k^\mathrm{Lan}\leq \frac{\log(n/\delta^2)}{4\sqrt{\epsilon}}
\]
iterations, with probability at least $1-\delta$. This is of course a highly conservative bound and in particular, the worst case matrices used to prove it vary with $k^\mathrm{Lan}$.

This means that computing one leading eigenpair of the matrix $X$ requires computing at most $k^\mathrm{Lan}$ matrix vector products (we can always restart the code in case of failure) plus $4nk^\mathrm{Lan}$ flops. When the matrix is dense, each matrix vector product costs $n^2$ flops, hence the total cost of computing one leading eigenpair of $X$ is
\[
O\left(\frac{n^2\log(n/\delta^2)}{4\sqrt{\epsilon}}\right)
\]
flops. When the matrix is sparse, the cost of each matrix vector product is $O(s)$ instead of $O(n^2)$, where $s$ is the number of nonzero coefficients in $X$. Idem when the matrix $X$ has rank $r<n$ and an explicit factorization is known (which is the case in the algorithms detailed in the previous section), in which case each matrix vector product costs $O(nr)$ which is the cost of two $n$ by $r$ matrix vector products, and the complexity of the Lanczos procedure decreases accordingly.

The numerical package ARPACK by \cite{Leho98} implements the Lanczos procedure with a reverse communication interface allowing the user to efficiently compute the matrix vector product $Xu_j$. However, it uses the implicitly shifted QR method instead of the more efficient MRRR algorithm to compute the Ritz pairs of the matrix $T_k\in\symm_k$.

\subsection{Other sampling techniques}


For completeness, we recall below another subsampling procedure in \cite{Achl07}. More recent ``volume sampling'' techniques produce improved error bounds (some with multiplicative error bounds) but the corresponding optimal sampling probabilities are much harder to compute, we refer the reader to \cite{Vemp09} for more details. The key idea behind this result is that, as the matrix dimension $n$ grows and given a fixed, scale invariant precision target $\|X\|_F/\epsilon$, the norm $\|X\|_\infty$ of individual coefficients in $X$ typically becomes negligible and we can randomly discard the majority of them while keeping important spectral features of $X$ mostly intact. 
\begin{lemma} \label{lem:rand-achl}
Given $X\in\symm_n$ and $\epsilon>0$, we define a subsampled matrix $S$ whose coefficients are independently distributed as:
\BEQ\label{eq:subsamp-achl}
S_{ij}=\left\{\BA{cl}
X_{ij}/p & \mbox{with probability $p$,}\\
0  & \mbox{otherwise.}\\
\EA\right.
\EEQ
when $i\geq j$, and $S_{ij}=S_{ji}$ otherwise. Assume that $1 \geq p\geq (8\log n)^4/n$, then
\[
\|X-S\|_2 \leq  4 \|X\|_\infty \sqrt{n/p}.
\]
with probability at least $1-\exp(-19(\log n)^4)$.
\end{lemma}
\begin{proof}
See \cite[Th.~1.4]{Achl07}.
\end{proof}

At first sight here, bounding the approximation error means letting the probability $p$  grow relatively fast as $n$ tends to infinity. However, because $\|X\|_\infty/\epsilon$ is typically much smaller than $\|X\|_F/\epsilon$, this subsampling ratio $p$ can often be controlled. Adaptive subsampling, i.e. letting $p$ vary with the magnitude of the coefficients in $X$, can further improve these results (see \cite[\S4]{Achl07} for details). The average number of nonzero coefficients in the subsampled matrix can be bounded using the structure of $X$. Note that the constants in this result are all very large (in particular,  $1 \geq p\geq (8\log n)^4/n$ implies $n\geq 10^9$) so despite its good empirical performance in low dimensions, the result presented above has to be understood in an asymptotic sense.

\section*{Acknowledgements}
The author is grateful to two anonymous referees for their numerous comments, and would like to acknowledge partial support from NSF grants SES-0835550 (CDI), CMMI-0844795 (CAREER), CMMI-0968842, a starting grant from the European Research Council (ERC project SIPA), a Peek junior faculty fellowship, a Howard B. Wentz Jr. award and a gift from Google.


\small{
\bibliographystyle{alpha}
\bibliography{Mainperso}}
\end{document}